\documentclass[11pt,titlepage]{article}
\usepackage{latexsym}
\usepackage{amsfonts}
\usepackage{amsmath}
\newcommand\blackslug{\hbox{\hskip 1pt \vrule width 4pt height 8pt depth 1.5pt
        \hskip 1pt}}
\newcommand\bbox{\hfill \quad \blackslug \medbreak}

\newtheorem{theorem}{}[section]

\newcommand{\Proof}{\noindent{\bf Proof.}\ \ }
\title{The Strong  EH-Property and the Erd\H{o}s-Hajnal Conjecture}
\author{
Krzysztof Choromanski
Columbia University\\
New York, NY, USA
}

\newcommand{\Keywords}{the Erd\H{o}s-Hajnal Conjecture, the regularity lemma, tournaments, EH-property, strong EH-property}

\begin{document}
\maketitle
\begin{abstract}

The Erd\H{o}s-Hajnal Conjecture states that for every $H$ there exists a constant $\epsilon(H)>0$ such that
every graph $G$ that does not contain $H$ as an induced subgraph contains a clique or a stable set of size at least
$|V(G)|^{\epsilon(H)}$. The Conjecture is still open. Some time ago its directed version was formulated (see:\cite{alon}). In the directed version graphs are replaced by tournaments, and cliques and stable sets by
transitive subtournaments. If the Conjecture is not true then the smallest counterexample is a prime tournament. For a long time the Conjecture was known only for finitely many prime tournaments. Recently in \cite{bcc} and \cite{choromanski2} the Conjecture was proven for the families of galaxies and constellations that contain infinitely many prime tournaments. In \cite{bcc} the Conjecture was also proven for all $5$-vertex tournaments.
We say that a tournament $H$ has the $EH$-property if it satisfies the Conjecture. In this paper we introduce the so-called \textit{strong EH-property} which enables us to prove the Conjecture for new prime tournaments, but what is even more interesting, provides a mechanism to combine tournaments satisfying the Conjecture to get bigger tournaments that do so and are not necessarily nonprime. We give several examples of families of tournaments constructed according to this procedure. The only procedure known before used to construct bigger tournaments satisfying the Conjecture from smaller tournaments satisfying the Conjecture was the so-called \textit{substitution procedure} (see: \cite{alon}). However an outcome of this procedure is always a nonprime tournament and, from what we have said before,  prime tournaments are those that play crucial role in the research on the Conjecture. Our method may be potentially used to prove the Conjecture for several new classes of tournaments.

\end{abstract}

\maketitle {\bf Keywords:} \Keywords

\section{Introduction}

We use $||$ to denote the size of the set. For a predicate $P$ we denote by $[P]$ the indicator function of $P$, i.e. $[P]=1$ if $P$ is true and is $0$ otherwise. By $(v_{1},...,v_{n})$ we denote an $n$-element vector with entries: $v_{1},...,v_{n}$. Let $G$ be a graph. We denote by $V(G)$ the set of its vertices. Sometimes instead of writing $|V(G)|$ we use shorter notation $|G|$. We call $|G|$ the \textit{size of G}. We denote by $E(G)$ the set of edges of a graph $G$. A \textit{clique} in the undirected graph is a set of pairwise adjacent vertices and an \textit{stable set} in the undirected graph is a set of pairwise nonadjacent vertices.
A \textit{tournament} is a directed graph  such that for every pair $v$ and $w$ of vertices, exactly one of the edges $(v,w)$ or $(w,v)$ exists. If $(v,w)$ is an edge of the tournament then we say that $v$ is $\textit{adjacent to}$ $w$ and $w$ is \textit{adjacent from} $v$. For two sets of vertices $V_{1}$, $V_{2}$ we say that $V_{1}$  is \textit{complete to} $V_{2}$ (or equivalently $V_{2}$ is \textit{complete from} $V_{1}$) if every vertex of $V_{1}$ is adjacent to every vertex of $V_{2}$.
A tournament is \textit{transitive} if it contains no directed cycle. For the set of vertices $V=\{v_{1},v_{2},...,v_{k}\}$ we say that an ordering $(v_{1},v_{2},...,v_{k})$ is \textit{transitive} if $v_{1}$ is adjacent to all other vertices of $V$, $v_{2}$ is adjacent to all other vertices of $V$ but $v_{1}$, etc. 
A tournament $H^{c}$ is a \textit{complement} of the tournament $H$ if $H^{c}$ is obtained from $H$ by reversing directions of all edges.
We denote by $C_{5}$ the unique tournament on $5$ vertices, where every vertex has two other vertices adjacent to it and two adjacent from it.
Let $\mathcal{H}$ be a family of tournaments.
If a tournament $T$ does not contain tournament $H$ as a subtournament for every $H \in \mathcal{H}$ then we say that $T$ is $\mathcal{H}$-\textit{free}. 
If $\mathcal{H}=\{H\}$ and $T$ is $\mathcal{H}$-free then we simply say that $T$ is $H$-\textit{free}.

A celebrated unresolved Conjecture of  Erd\H{o}s and Hajnal states that:  

\begin{theorem}
\label{EHConun}
For every undirected graph $H$ there exists $\epsilon(H)>0$ such that every $n$-vertex undirected graph that does not contain $H$ as an induced subgraph contains a clique or a stable of size at least $n^{\epsilon(H)}$.
\end{theorem}

In 2001 Alon, Pach and Solymosi proved (see: \cite{alon}) that Conjecture~\ref{EHConun} has an equivalent directed version, where undirected graphs are replaced by tournaments and cliques and stable sets by transitive subtournaments.

The equivalent directed version (see: \cite{alon}) states that:

\begin{theorem}
\label{EHCon}
For every tournament $H$ there exists $\epsilon(H)>0$ such that every $n$-vertex $H$-free tournament contains a transitive subtournament of size at least $n^{\epsilon(H)}$.
\end{theorem}

If for a graph $H$ there exists $\epsilon(H)>0$ then we say that \textit{$H$ satisfies the Erd\H{o}s-Hajnal Conjecture} or equivalently: \textit{has EH-property}. For a given family of tournaments $\mathcal{H}$, if there exists $\epsilon(\mathcal{H})>0$ such that every $\mathcal{H}$-free tournament $T$ contains a transitive subtournament of size at least $|T|^{\epsilon(\mathcal{H})}$ then we say that \textit{$\mathcal{H}$ has EH-property}. 
Obviously, the analogous definition may be stated for the family $\mathcal{H}$ of undirected graphs. However, since in this paper we focus on the directed version of the Conjecture, we do not need the corresponding definition in the undirected setting.

Whenever we work with tournaments it is convenient to fix some ordering of its vertices. Below we give some definitions introducing that approach.
We also remind the definition of the family of galaxies that may be also found in \cite{bcc}.

Let $T$ be a tournament, and let $(v_1, \ldots, v_{|T|})$ be an ordering of its 
vertices; denote this ordering by $\theta$.  We say that an edge $(v_{j},v_{i})$ of $T$ is a \textit{backward edge} under this ordering if $i<j$.  
The \textit{graph of backward edges}  under this ordering,
denoted by $B(T, \theta)$, has vertex set $V(T)$,  and  
$v_i v_j \in E(B(T, \theta))$ if and only if  $(v_i,v_j)$ or 
$(v_j,v_i)$ is a backward edge of $T$ under the ordering $\theta$. 
For an integer $t$, we call the graph $K_{1,t}$ a {\em star}. Let $S$ be a
star with vertex set $\{c, l_1, \ldots, l_t\}$, where $c$ is adjacent to $l_1, \ldots, l_t$. We call $c$ the {\em center of the star}, and
$l_1, \ldots, l_t$ {\em the leaves of the star}. 
Note that in the case $t=1$ we may choose arbitrarily any one of the two vertices to be the center of the star, and the other vertex is then considered to be the leaf. 

A {\em right star}
in $B(T, \theta)$ is an induced subgraph with vertex set 
$\{v_{i_0}, \ldots, v_{i_t}\}$, such that \\
$B(T,\theta)|\{v_{i_0}, \ldots, v_{i_t}\}$ is a star with center $v_{i_t}$, 
and  $i_t > i_0, \ldots, i_{t-1}$.  In this case we also  
say that $\{v_{i_0}, \ldots, v_{i_t}\}$ is  a right star in $T$.
A {\em left star}
in $B(T, \theta)$ is an induced subgraph with vertex set 
$\{v_{i_0}, \ldots, v_{i_t}\}$, such that 
$B(T,\theta)|\{v_{i_0}, \ldots, v_{i_t}\}$ is a star with center $v_{i_0}$, 
and  $i_0 < i_1, \ldots, i_t$.    In this case we also  
say that $\{v_{i_0}, \ldots, v_{i_t}\}$ is a  left star in $T$.
Finally, a {\em star}   in $B(T,\theta)$,  is a left star or a right star.

A tournament $T$ is a \textit{galaxy} if there exists an ordering $\theta$ of its vertices such that every connected component of $B(T, \theta)$ is  either a star or a 
singleton, and
\begin{itemize}
\item 
no center of a star appears in the ordering between two leaves of another star.
\end{itemize}
We call such an ordering a {\em galaxy ordering} of $T$. Let 
$\Sigma_1, \ldots, \Sigma_l$ be the non-singleton components of 
$B(T, \theta)$. We say that $\Sigma_1, \ldots, \Sigma_l$ 
are the {\em stars of $T$ under theta}. If $V(T)=\bigcup_{i=1}^l V(\Sigma_l)$, 
we say that $T$ is a {\em regular} galaxy.

For a tournament $H$ a subset $S \subseteq V(H)$ is called \textit{homogeneous} if for every $v \in V(H) \backslash S$ the following holds: either $\forall_{w \in S} (w,v)$ is an edge or $\forall_{w \in S} (v,w)$ is an edge. A homogeneous set $S$ is called \textit{nontrivial} if $|S|>1$ and $S \neq V(H)$. A tournament is called \textit{prime} if it does not have nontrivial homogeneous sets. A tournament is called \textit{nonprime} if it is not prime.
Analogous definitions may be introduced for undirected graphs but once more we will not do it since we will not use them later in the paper.
 
The following theorem is an immediate corollary of the results given in \cite{alon}. Applied for tournaments, shows why prime tournaments are important.

\begin{theorem}
If Conjecture~\ref{EHCon} is false then the smallest counterexample is prime.
\end{theorem}

The theorem follows from the \textit{substitution procedure} presented in \cite{alon} and adapted to the directed setting. The procedure takes as an input two tournaments satisfying the Conjecture and outputs a bigger tournament that satisfies it too. For a long time the procedure was the only method to construct infinitely many tournaments satisfying the Conjecture. However it can be easily noticed that an outcome of the procedure is a tournament that always has nontrivial homogeneous sets. Therefore an outcome is always nonprime. The question arises whether it is possible to show that there are infinitely many prime graphs satisfying the Conjecture. In the undirected case it is still open since so far the Conjecture is known only for some undirected prime graphs of at most $5$ vertices. In the directed setting very recently the Conjecture was proven for new families of tournaments - galaxies and constellations that contain infinitely many prime tournaments. The Conjecture was also proven for all tournaments on at most $5$ vertices (see: \cite{bcc} and \cite{choromanski2}). 

In this paper we introduce new property of a family of tournaments $\mathcal{H}$ called the \textit{strong EH-property}. It will enable us to prove the Conjecture for new families of tournaments that are not contained in the familes of galaxies and constellations. Those methods can be also used to prove the Conjecture for galaxies, constellations, all tournaments on at most $5$ vertices and all but one tournament on $6$ vertices. Furthermore, presented techniques enable us to give new methods of combining tournaments satisfying the Conjecture to get bigger tournaments which also do so. They may be useful in proving the Conjecture for many new families of tournaments. To the best of our knowledge, our method is the only one apart from the substitution procedure that produce infinitely many tournaments with the EH-property. But in contrast to the substitution procedure, those tournaments are not necessarily nonprime. Therefore methods presented by us may be used to produce infinitely many prime tournaments with EH-property that are neither galaxies nor constellations.

This paper is organized as follows:
\begin{itemize}
\item in Section 2 we introduce definition of the strong EH-property and several technical definitions used by us later in the paper,
\item in Section 3 we prove theorems about the strong EH-property,
\item in Section 4 we give several applications of the strong EH-property and theorems presented in the previous section, proving the Conjecture for new families of tournaments,
\item in Section 5 we give some final remarks.
\end{itemize}

\section{Preliminaries}

In this section we give the definition of the strong EH-property. However, before doing it, we need to introduce several technical definitions.

\subsection{m-sequences}

For a $\{0,1\}$-vector $v$ denote $\zeta^{v}(i)=|{j:j \leq i, v(j)=1}|$.
Denote by $tr(T)$ the size of the biggest transitive subtournament of $T$.

Let $X,Y \subseteq V(T)$ be disjoint. Denote by $e_{X,Y}$ the number of directed edges $(x,y)$, where $x \in X$ and $y \in Y$.
The \textit{directed density from X to Y} is defined as 
$d(X,Y)=\frac{e_{X,Y}}{|X||Y|}.$ 

Let $T$ be a tournament and let $S_{1},S_{2},...,S_{k} \subseteq V(T)$ s.t. $S_{i} \bigcap S_{j} = \emptyset$ for $1 \leq i < j \leq k$. Let $v$ be a $\{0,1\}$-vector of length $k$ and let $\eta \in N^{\zeta^{v}(k)}$. For a given $c>0$, $0< \lambda < 1$ we say that a sequence $(S_{1},...,S_{k})$ is the $(v,\eta,c,\lambda)$-m-sequence
if the following holds:
\begin{itemize}
\item if $v_{i}=1$ then $S_{i}$ induces a transitive subtournament,
\item if $v_{i}=0$ then $|S_{i}| \geq c|T|$,
\item each $S_{i}$ for which $v_{i}=1$ is partitioned into $\eta(\zeta^{v}(i))$ subsets: $S_{i,1},...,S_{i,\eta(\zeta^{v}(i))}$ s.t. $S_{i,j_{1}}$ is complete to $S_{i,j_{2}}$ for $1 \leq j_{1} < j_{2} \leq \eta(\zeta^{v}(i))$ and $|S_{i,j}| \geq c tr(T)$ for $j=1,2,...,\eta(\zeta^{v}(i))$,
\item $d(S_{i},S_{j}) \geq 1 - \lambda$ for $1 \leq i < j \leq k$. 
\end{itemize}

If we do not care about parameters, we simply say: $m$-sequence. We say that a $m$-sequence $\chi$ is \textit{strong} if:
\begin{itemize}
\item $d(\{v\},S_{j}) \geq 1 - \lambda$ for every $i<j$, $v \in S_{i}$,
\item $d(S_{i},\{v\}) \geq 1 - \lambda$ for every $j>i$, $v \in S_{j}$.
\end{itemize}

For the vectors $v,\eta$ as in the definition of the $m$-sequence we denote 

$$\gamma^{v,\eta}(i) = |\{j \leq i: v(j)=0\}| + \sum_{j \leq i, v(j)=1} \eta(\zeta^{v}(j)),$$ for $i=1,2,...,|v|$.

Whenever we will talk about $m$-sequences, it will always be in the context of the tournament $T$ from which vertices of the $m$-sequence are taken. Thus by the subtournament of the $m$-sequence $\chi=(S_{1},...,S_{l})$ we mean the subtournament of the tournament $T$ induced by $\bigcup_{i=1}^{l}S_{i}$.
We denote $V(\chi)=\bigcup_{i=1}^{l}S_{i}$. 
We will also use the longer representation of the $m$-sequence $\chi=(S_{1},...,S_{l})$, denoted by $l(\chi)$. The sequence $l(\chi)$ is obtained from $\chi$ by replacing each $S_{i}$ s.t. $v_{i}=1$ by a sequence: $S_{i,1},...,S_{i,\eta(\zeta^{v}(i))}$. Thus the length of $l(\chi)$ is $|\{j:v_{j}=0\}| + \sum_{i} \eta(i)$.

For vectors $v_{1},v_{2} \in \{0,1\}^{m}$, $\eta_{1} \in N^{\zeta^{v_{1}}(m)}$, $\eta_{2} \in N^{\zeta^{v_{2}}(m)}$  we say that $(v_{1},\eta_{1}) \subseteq (v_{2},\eta_{2})$ if:

\begin{itemize}
\item $v_{1}=(v_{2}(i_{1}),...,v_{2}(i_{t}))$ for some $i_{1} < i_{2} < ... < i_{t}$ , i.e. $v_{1}$ is a subsequence of $v_{2}$,
\item if $v_{1}(j)=1$ then $\eta_{1}(\zeta^{v_{1}}(j)) \leq \eta_{2}(\zeta^{v_{2}}(i_{j}))$.
\end{itemize}

Let $T$ be a tournament and let $A,B \in V(T)$. Let $A \cap B = \emptyset$. 
Assume that $|A| \geq c|T|$ for some $c>0$. Assume furthermore that $|B| \geq c|T|$ or $B$ induces a transitive subtournament and $|B| \geq c tr(T)$. Assume that  $d(A,B)=1$ or $d(B,A)=1$.
Then we say that $(A,B)$ is a \textit{$c$-strong pair}. If a parameter $c$ is not important then we simply say that $(A,B)$ is a \textit{strong pair}.

We will introduce now few more definitions, among them the crucial definition of the strong EH-property. 

Let $v,\eta$ be vectors from the definition of the $m$-sequence. Let $\mathcal{H} = \{H_{1},...,H_{t}\}$ be a finite family of tournaments.
Denote $V(H_{i})=\{h^{i}_{1},...,h^{i}_{|H_{i}|}\}$.
Let $k = \sum_{i} \eta(i) + \sum_{j} [v(j)=0]$. Let $\mathcal{F}=\{f_{1},...,f_{t}\}$, where for $i=1,...,t$ $f_{i}:V(H_{i}) \rightarrow \{1,...,k\}$ is an injective function. 
We say that $\mathcal{H}$ \textit{has $(\mathcal{F},v,\eta)$-strong EH-property} (or just: \textit{has strong EH-property} when parameters are not important) if: \\

$\forall_{c_{1}>0}$ $\exists \lambda_{0}(c_{1})>0, c_{2}(c_{1})>0$ s.t.
if $\chi$ is a $(v,\eta,c_{1},\lambda)$-m-sequence for $0 \leq \lambda \leq \lambda_{0}(c_{1})$ then the following holds:

\begin{itemize}
\item $\exists i, v_{1},...,v_{|H_{i}|}$ s.t. $\{v_{1},...,v_{|H_{i}|}\}$ induces a subtournament of $\chi$ isomorphic to $H_{i}$, where: $v_{j} \in T_{f_{i}(h^{i}_{j})}$ for $j=1,2,...,|H_{i}|$, $l(\chi) = (T_{1},...,T_{k})$
and the isomorphism is defined by the mapping $v_{i} \rightarrow h^{i}_{j}$ or 
\item $\chi$ contains a $c_{2}$-strong pair.
\end{itemize}

Let $T$ be a tournament and let $\chi$ be a $m$-sequence in $T$.
We say that a tournament $H$ with $V(H)=\{h_{1},...,h_{|H|}\}$ is $(f,\lambda)$-well-embedded in $\chi$, where
$f: V(H) \rightarrow \{1,...,k\}$ and $l(\chi) = (T_{1},...,T_{k})$ if:

\begin{itemize}
\item $\exists v_{1},...v_{|H|}$ s.t. $v_{i} \in T_{f(h_{i})}$,  $\{v_{1},...,v_{|H_{i}|}\}$ induces a subtournament of $\chi$ isomorphic to $H$, the isomorphism is defined by the mapping $v_{i} \rightarrow h_{i}$ and
\item for every $i$ s.t. $v_{1},...,v_{|H|} \notin T_{i}$ and $j \in \{1,...,|H|\}$ we have: $d(\{v_{j}\},T_{i}) \geq 1-\lambda$ if $f(h_{j})<i$
and $d(T_{i},\{v_{j}\}) \geq 1-\lambda$ if $f(h_{j})>i$.
\end{itemize}

Let us take denotations from the definition of the strong EH-property.

We say that \textit{$\mathcal{H}$ has $(\mathcal{F},v,\eta)$-super-strong EH-property} if: \\

$\forall_{c_{1}>0}$ $\exists \lambda_{0}(c_{1})>0, c_{2}(c_{1})>0$ s.t.
if $\chi$ is a $(v,\eta,c_{1},\lambda)$-m-sequence for $0 \leq \lambda \leq \lambda_{0}(c_{1})$ then the following holds:

\begin{itemize}
\item $\exists i$ s.t. $H_{i}$ is $(f_{i},\frac{1}{2k})$-well-embedded in $\chi$ or
\item $\chi$ contains a $c_{2}$-strong pair.
\end{itemize}

Let $v_{1},v_{2},\eta_{1},\eta_{2}$ be vectors as in the definition of the $m$-sequence. Assume that $(v_{1},\eta_{1}) \subseteq (v_{2},\eta_{2})$.
Let  $k_{1} = \sum_{i} \eta_{1}(i) + \sum_{j} [v_{1}(j)=0]$ and $k_{2} = \sum_{i} \eta_{2}(i) + \sum_{j} [v_{2}(j)=0]$. 
Denote $v_{1}=(v_{2}(i_{1}),...,v_{2}(i_{t}))$ for some $i_{1} < ... < i_{t}$.
Let $s \in \{1,2,...,k_{2}\}$. Let $j^{s}$ be the smallest index $j$ such that $\gamma^{v_{2},\eta_{2}}(j) \geq s$. Assume that $j^{s}=i_{r}$ for some $1 \leq r \leq t$. Let $\beta(s)=\{w \in \{1,...,k_{1}\} |\gamma^{v_{1},\eta_{1}}(w)=r\}$.
Let $\mathcal{H}=\{H_{1},...,H_{r}\}$ be a finite family of tournaments.
Denote $V(H_{i})=\{h^{i}_{1},...,h^{i}_{|H_{i}|}\}$ for $i=1,...,r$ and
let $\mathcal{F}_{1}=\{f^{1}_{1},...,f^{1}_{r}\}$ be a family of functions such that $f^{1}_{i}:V(H_{i}) \rightarrow \{1,...,k_{1}\}$. Let $\mathcal{F}_{2}=\{f^{2}_{1},...,f^{2}_{r}\}$ be a family of functions such that $f^{2}_{i}:V(H_{i}) \rightarrow \{1,...,k_{2}\}$ and besides:

\begin{itemize}
\item $\beta(f^{2}_{i}(h^{i}_{j}))$ is defined for every $i=1,...,r$, $j=1,...,|H_{i}|$,
\item  $f^{1}_{i}(h^{i}_{j}) \in \beta(f^{2}_{i}(h^{i}_{j}))$ for every $i=1,...,r$, $j=1,...,|H_{i}|$,
\item $f^{1}_{i_{1}}(h^{i_{1}}_{j_{1}}) \leq f^{1}_{i_{2}}(h^{i_{2}}_{j_{2}})$ iff $f^{2}_{i_{1}}(h^{i_{1}}_{j_{1}}) \leq f^{2}_{i_{2}}(h^{i_{2}}_{j_{2}})$.  
\end{itemize}

Then we say that triple \textit{$(\mathcal{F}_{2},v_{2},\eta_{2})$ EH-extends $(\mathcal{F}_{1},v_{1},\eta_{1})$}. It is easy to see that if $\mathcal{H}$ has $(\mathcal{F}_{1},v_{1},\eta_{1})$-strong EH-property then it also has $(\mathcal{F}_{2},v_{2},\eta_{2})$-strong EH-property. 

\begin{theorem}
\label{EHext_theorem}
If $\mathcal{H}$ has $(\mathcal{F}_{1},v_{1},\eta_{1})$-strong EH-property then it also has $(\mathcal{F}_{2},v_{2},\eta_{2})$-strong EH-property, where  $(\mathcal{F}_{2},v_{2},\eta_{2})$ EH-extends $(\mathcal{F}_{1},v_{1},\eta_{1})$.  
\end{theorem}

\Proof
Take a $(v_{2},\eta_{2},c_{2},\lambda)$-$m$-sequence $\chi_{2}=(S^{2}_{1},...,S^{2}_{m})$ for some arbitrary $c_{2}>0$, where: $\lambda \leq c_{2}^{2}\lambda_{0}(c_{2})$ and $\lambda_{0}(c_{2})$ is as in the definition of the $(\mathcal{F}_{1},v_{1},\eta_{1})$-strong EH-property. Let $l(\chi_{2})=(T_{1},...,T_{k})$. Note that directly from the definition of the EH-extension we know that there exists a $(v_{1},\eta_{1},c_{2},\frac{\lambda}{c_{2}^{2}})$-$m$-sequence $\chi_{1}=(S^{1}_{i_{1}},...,S^{1}_{i_{t}})$ with $l(\chi_{2})=(T_{j_{1}},...,T_{j_{r}})$ for some $1 \leq i_{1} < ... < i_{t} \leq m$, $1 \leq j_{1} < ... < j_{r} \leq k$, where $S^{1}_{i_{k}} \subseteq S^{2}_{i_{k}}$ for $k=1,...,t$ and such that:
\begin{itemize}
\item an embedding of the tournament $H \in \mathcal{H}$ in $\chi_{2}$ from the definition of the $(\mathcal{F}_{2},v_{2},\eta_{2})$-strong EH-property corresponds to the embedding of the tournament $H \in \mathcal{H}$ in $\chi_{1}$ from the definition of the $(\mathcal{F}_{1},v_{1},\eta_{1})$-strong EH-property. 
\end{itemize}
It sufficies to note that for chosen $\lambda$ we can use $(\mathcal{F}_{1},v_{1},\eta_{1})$-strong EH-property and either get a desired embedding or a strong pair. That completes the proof.
\bbox

\section{Strong EH-property and the Erd\H{o}s-Hajnal Conjecture}

In this section we prove several results connecting strong EH-property with the Erd\H{o}s-Hajnal Conjecture.

We start with the very useful technical lemma:
 
\begin{theorem}
\label{technical_lemma}
Let $T$ be a tournament and let $X,Y \subseteq V(T)$, $X \bigcap Y = \emptyset$.
Assume that $d(X,Y) \geq 1 - \lambda$. Let $X_{1} \subseteq X$, $Y_{1} \subseteq Y$. Assume that $|X_{1}| \geq c_{1}|X|$, $|Y_{1}| \geq c_{2}|Y|$.
Then $d(X_{1},Y_{1}) \geq 1 - \frac{\lambda}{c_{1}c_{2}}$.
\end{theorem} 

\Proof
The number of directed edges $n^{X,Y}_{e}$ going from $Y$ to $X$ is $(1-d(X,Y))|X||Y|$. Thus $n^{X,Y}_{e} \leq \lambda |X||Y|$. The number of directed edges $n^{X_{1},Y_{1}}_{e}$ going from $Y_{1}$ to $X_{1}$ is $(1-d(X_{1},Y_{1}))|X_{1}||Y_{1}|$. Since  $n^{X_{1},Y_{1}}_{e} \leq n^{X,Y}_{e} $, the result follows.  
\bbox

\subsection{Strong EH-property implies the Erd\H{o}s-Hajnal Conjecture}

We show now that strong EH-property implies EH-property.

\begin{theorem}
\label{EHmain}
Let $\mathcal{H}$ be a finite nonempty family of tournaments. If $\mathcal{H}$ has $(\mathcal{F},v,\eta)$-strong EH-property then $\mathcal{H}$ has EH-property.
\end{theorem}

\Proof
Let $T$ be a $\mathcal{H}$-free $n$-vertex tournament. The proof is by induction on $n$. 
Let $N(v,\eta)$ be the smallest integer such that there exists $c_{1}>0$ with the following property: every $n$-vertex $\mathcal{H}$-free tournament $T$, where $n \geq N(v,\eta)$, contains a $(v,\eta,c_{1},\lambda)$-$m$-sequence $\chi$ for some $0 \leq \lambda \leq \lambda_{0}(c_{1})$ ($\lambda_{0}(c_{1})$ as in the definition of the strong EH-property). It can be proven (standard application of the regularity lemma for tournaments, see: \cite{shapira}) that $N(v,\eta)$ is finite. Formal proof may be found in \cite{bcc}, therefore we will not repeat it now. We will prove that every $\mathcal{H}$-free tournament $T$ contains a transitive subtournament of size at least $|T|^{\epsilon(N,c_{1})}$, 
for: $\epsilon=\epsilon(N,c_{1}) = \min(\log_{2}(N),\log_{c_{2}(c_{1})}(1-c_{2}(c_{1})),\log_{c_{2}}(\frac{1}{2}))$, where $c_{2}(c_{1})$ is as in the definition of the strong EH-property.
Note that every tournament $T^{'}$ of at most $N$ vertices contains a transitive subtournament of size at least $|T^{'}|^{\epsilon}$, since $\epsilon \leq \log_{2}(N)$.
Now assume that every tournament $T^{'}$ on at most $n_{0}$ vertices contains a transitive subtournament of size at least $|T^{'}|^{\epsilon}$.
We may thus assume that $T$ has $n=n_{0}+1$ vertices. Since $T$ is $\mathcal{H}$-free, we can conclude, by the previous remark, that $T$ contains a $(v,\eta,c_{1},\lambda)$-$m$-sequence $\chi$. Now, since $T$ is $\mathcal{H}$-free and has strong EH-property, we notice that $T$ contains a $c_{2}(c_{1})$-strong pair $(A,B)$.
Assume first that $B$ induces a transitive subtournament.
Then, by induction we have: $tr(A) \geq |A|^{\epsilon}$. Since: $tr(T) \geq tr(A) + |B|$, from the definition of the $c_{2}$-strong pair we get:
$tr(T) \geq (c_{2}n)^{\epsilon} +c_{2} tr(T)$. Thus we have: $tr(T) \geq \frac{c_{2}^{\epsilon}n^{\epsilon}}{1-c_{2}}$. Thus, since $\epsilon \leq \log_{c_{2}}(1-c_{2})$, we get: $tr(T) \geq n^{\epsilon}$ and we are done.
Assume now that $|B| \geq c_{2}n$. As, in the previous case, we have by induction: $tr(A) \geq |A|^{\epsilon}$. We also have: $tr(B) \geq |B|^{\epsilon}$.  Since: $tr(T) \geq tr(A) + tr(B)$, from the definition of the $c_{2}$-strong pair we get:
$tr(T) \geq 2(c_{2}n)^{\epsilon}$. Thus, since: $\epsilon \leq \log_{c_{2}}(\frac{1}{2})$, we get: $tr(T) \geq n^{\epsilon}$ and that  completes the proof.
\bbox

From Theorem~\ref{EHmain} we immediately get the following corollary:

\begin{theorem}
\label{EHStrongConj}
Let $\mathcal{H}=\{H_{1},...,H_{m}\}$, where all $H_{i}^{s}$ contain some tournament $H$. Assume furthermore that $\mathcal{H}$ has $(\mathcal{F},v,\eta)$-strong EH-property. Then $\mathcal{H}$ has EH-property thus $H$ satisfies the Erd\H{o}s-Hajnal Conjecture. 
\end{theorem}

We will prove now that strong EH-property is in fact equivalent to the super-strong EH-property. This fact will turn out to be very useful later when we will consider so-called \textit{product tournaments}.

\begin{theorem}
\label{strong_superstrong}
Let $\mathcal{H}$ be a finite family of tournaments. If $\mathcal{H}$ has $(\mathcal{F},v,\eta)$-strong EH-property then it has $(\mathcal{F},v,\eta)$-super-strong EH-property.
\end{theorem}

\Proof
Let $l(\chi)=(T_{1},...,T_{m})$.
Take an arbitrary constant $c_{1}>0$ and consider $(v,\eta,c_{1},\lambda)$-m-sequence $\chi$, where $\lambda = \min(\frac{c_{1}^{2}\lambda_{0}(c_{1})}{4m},\frac{1}{4km})$, where $\lambda_{0}(c_{1})$ is as in the definition of the strong EH-property.
Note that from Theorem~\ref{technical_lemma} we know that $d(T_{i},T_{j}) \geq 1 - \frac{\lambda}{c_{1}^{2}}$ for $i < j$. Let $\hat{\lambda}=\frac{\lambda}{c_{1}^{2}}$. Let $C=2m$.
Let $S^{i}_{j}$ be a subset of $T_{i}$ consisting of those vertices $v$ in $T_{i}$ that satisfy:
\begin{itemize}
\item $d(\{v\},T_{j}) \geq 1 - C\hat{\lambda}$ if $i < j$ or
\item $d(T_{j},\{v\}) \geq 1 - C\hat{\lambda}$ if $j < i$.
\end{itemize}

From Theorem~\ref{technical_lemma} we get: $|S^{i}_{j}| \geq (1-\frac{1}{C})|T_{i}|$. Denote $S^{i} = \bigcap_{j \neq i} S^{i}_{j}$.
We have: $|S^{i}| \geq (1-\frac{m}{C})|T_{i}|$. Replacing in the sequence $(T_{1},...,T_{m})$ every $T_{i}$ by $S^{i}$ we form a new $m$-sequence $\chi^{'}$ which is a $(v,\eta,c^{'}_{1},\lambda^{'})$-strong-m-sequence, where:
\begin{itemize}
\item $c^{'}_{1} = c_{1}(1-\frac{m}{C})$,
\item $\lambda^{'} = \frac{C \lambda^{'}}{1-\frac{m}{C}}=\frac{C\frac{\lambda}{c_{1}^{2}}}{1-\frac{m}{C}}$. 
\end{itemize}

Now we use the fact that $\mathcal{H}$ has $(\mathcal{F},v,\eta)$-strong EH-property. So, since $\lambda^{'} \leq \lambda_{0}(c_{1})$, then we have either:

\begin{itemize}
\item $\chi^{'}$ contains $c_{2}$-strong pair for some $c_{2} > 0$ or
\item $\exists i, v_{1},...,v_{|H_{i}|}$ s.t. $\{v_{1},...,v_{|H_{i}|}\}$ induces a subtournament of $\chi^{'}$ isomorphic to $H_{i}$ for some $H_{i} \in \mathcal{H}$, where: $V(H_{i})=\{h^{i}_{1},...,h^{i}_{|H_{i}|}\}$, $v_{j} \in T_{f_{i}(h^{i}_{j})}$ for $j=1,2,...,|H_{i}|$, 
and the isomorphism is defined by the mapping $v_{i} \rightarrow h^{i}_{j}$.
\end{itemize}

In the first case we are obviously done. Now assume the second case. Since $\chi^{'}$ is a strong $m$-sequence, then we see that the embedding defined by $\{v_{1},...,v_{|H_{i}|}\}$ is in fact a $(f_{i},C\lambda)$-well-embedding. Thus, since $C\lambda \leq \frac{1}{2k}$, we are done. 
\bbox

\subsection{Product tournaments}

In this subsection we will describe general method that enables us to construct prime tournaments satisfying the Conjecture from smaller tournaments that satisfy the Conjecture, by combining them in some specific way. \\

Let $H_{1},H_{2}$ be two tournaments. Let us consider two injective functions $f_{1}:V(H_{1}) \rightarrow \mathbb{N}, f_{2}:V(H_{2}) \rightarrow \mathbb{N}$.
Assume furthermore that $\forall_{h^{1} \in V(H_{1}),h^{2} \in V(H_{2})}$ we have: $f_{1}(h^{1}) \neq f_{2}(h^{2})$. We shortly denote this last condition by: $<f_{1},f_{2}>=0$ (For two families of functions: $\mathcal{F}$ and $\mathcal{G}$ we say that $<\mathcal{F},\mathcal{G}>=0$ if $\forall_{f \in \mathcal{F}, g\in \mathcal{G}} <f,g>=0$). Denote by $\theta_{1}$ the ordering of the vertices of $V(H_{1})$ induced by an increasing values of $f_{1}$ on  $V(H_{1})$ and by $\theta_{2}$ the ordering of the vertices of $V(H_{2})$ induced by an increasing values of $f_{2}$ on $V(H_{2})$.
Now let us define \textit{the product $H$ of $H_{1}$ and $H_{2}$ under orderings $\theta_{1}$ and $\theta_{2}$} as follows:
\begin{itemize}
\item $V(H)=V(H_{1}) \cup V(H_{2})$,
\item under ordering $\theta$ of $V(H)$ induced by $f_{1},f_{2}$, where: $<f_{1},f_{2}>=0$, the backward edges of $H$ are exactly the backward edges of $H_{1}$ under $\theta_{1}$ and the backward edges of $H_{2}$ under $\theta_{2}$.
\end{itemize}
We denote this product tournament $H$ by $H^{f_{1}}_{1} \oplus H^{f_{2}}_{2}$.

We will also need a notion of the product of families of tournaments. Before defining it we will give few more useful definitions. Let $H_{1},H_{2}$ be two tournaments with disjoint vertex sets and let $f:V(H_{1}) \rightarrow \mathbb{N},g: V(H_{2}) \rightarrow \mathbb{N}$ be two functions. 
Let us define on $V(H_{1}) \cup V(H_{2})$ \textit{the product of $f$ and $g$}:
$$ 
(f \oplus g)(h)=\left\{ 
\begin{array}{l l}
   f(h) & \quad \text{if $h \in V(H_{1})$}\\
   g(h) & \quad \text{if $h \in V(H_{2})$.}
\end{array}
\right.
$$

Let $\mathcal{F},\mathcal{G}$ be two families of functions defined on vertex-sets of tournaments such that $\forall_{f \in \mathcal{F},g \in \mathcal{G}}$ $f$ and $g$ are defined on disjoint domains. We define \textit{the product of $\mathcal{F}$ and $\mathcal{G}$} as follows: 
$$\mathcal{F} \oplus \mathcal{G} = \{f \oplus g: f \in \mathcal{F}, g \in \mathcal{G}\}.$$
Later on whenever we will use the product of the families of tournaments we will always assume the setting given above.

Take again two families of functions $\mathcal{F}=\{f_{1},...,f_{r}\},\mathcal{G}=\{g_{1},...,g_{s}\}$ as above.
Assume that $\mathcal{H}_{1}=\{H_{1,1},...,H_{1,r}\}$ has $(\mathcal{F},v,\eta)$-strong EH-property and that $\mathcal{H}_{2}=\{H_{2,1},...,H_{2,s}\}$ has $(\mathcal{G},v,\eta)$-strong EH-property for some vectors: $v,\eta$. Then we define \textit{the product $\mathcal{H}_{1}^{\mathcal{F}} \oplus \mathcal{H}_{2}^{\mathcal{G}}$  of families of tournaments $\mathcal{H}_{1}$ and $\mathcal{H}_{2}$ under $\mathcal{F}$ and $\mathcal{G}$} as follows: 

$$\mathcal{H}_{1}^{\mathcal{F}} \oplus \mathcal{H}_{2}^{\mathcal{G}} = \{H_{1,i}^{f_{i}} \oplus H_{2,j}^{g_{j}} | H_{1,i} \in \mathcal{H}_{1}, H_{2,j} \in \mathcal{H}_{2}\}.$$

Product tournaments are especially interesting in the context of the Conjecture because of the following theorem:

\begin{theorem}
\label{product_theorem}
Assume that a finite family of tournaments $\mathcal{H}_{1}$ has $(\mathcal{F},v,\eta)$-strong EH-property and that a finite family of tournaments $\mathcal{H}_{2}$ has $(\mathcal{G},v,\eta)$-strong EH-property.
Assume furthermore that $<\mathcal{F},\mathcal{G}>=0$. Then a family $\mathcal{H}=\mathcal{H}^{\mathcal{F}}_{1} \oplus \mathcal{H}^{\mathcal{G}}_{2}$ has the
$(\mathcal{F} \oplus \mathcal{G},v,\eta)$-strong EH-property. 
\end{theorem}

Before proving this theorem we briefly discuss some of its consequences.

Note first that from Theorem~\ref{EHStrongConj} and Theorem~\ref{product_theorem} we immediately have the following corollary:

\begin{theorem}
If $\mathcal{H}_{1}$ has $(\mathcal{F},v,\eta)$-strong EH-property and  $\mathcal{H}_{2}$ has $(\mathcal{G},v,\eta)$-strong EH-property then $\mathcal{H}_{1}^{\mathcal{F}} \oplus \mathcal{H}_{2}^{\mathcal{G}}$ satisfies the Erd\H{o}s-Hajnal Conjecture.
\end{theorem}

From Theorem~\ref{product_theorem} and the fact that EH-extensions preserve strong EH-property we also easily obtain the following result:

\begin{theorem}
\label{extension_theorem}
Let $\Omega=\{\mathcal{H}_{1},\mathcal{H}_{2},...\}$ be a family of families of tournaments. Assume that $\mathcal{H}_{i}$ has $(\mathcal{F}^{\mathcal{H}_{i}},v^{\mathcal{H}_{i}},\eta^{\mathcal{H}_{i}})$-strong EH-property for $i=1,2,...$. Denote by $\hat{\Omega}$ the closure of $\Omega$ obtained under taking EH-extensions and applying $\oplus$ operation to pairs of families of tournaments. Then every family of $\hat{\Omega}$ satisfies the Erd\H{o}s-Hajnal Conjecture. 
\end{theorem}

We will now introduce an important class of tournaments called \textit{regular tournaments} that satisfy the Erd\H{o}s-Hajnal Conjecture. This class contains all prime tournaments for which the Conjecture has been proven so far, such as: prime galaxies (see:\cite{bcc}) and prime constellations (see: \cite{choromanski2}). But this class is much larger. In fact it does not give us a well-defined new family of tournaments satisfying the Conjecture. It rather gives a new mechanism that enables us to prove the Conjecture for new families of prime tournaments by combining in a very specific way other tournaments with the EH-property. \\

Tournament $H$ is \textit{regular} if $\exists_{\mathcal{H}_{1},\mathcal{H}_{2},\mathcal{F},\mathcal{G}}$ (we may have: $\mathcal{H}_{1}=\emptyset$ or $\mathcal{H}_{2}=\emptyset$) s.t. $\mathcal{H}_{1}^{\mathcal{F}} \oplus \mathcal{H}_{2}^{\mathcal{G}} = \{H_{1},...,H_{|\mathcal{H}_{1}^{\mathcal{F}} \oplus \mathcal{H}_{2}^{\mathcal{G}}|}\}$, where each $H_{i}$ contains $H$ as a subtournament and besides: $\mathcal{H}_{1}$ has $(\mathcal{F},v,\eta)$-strong EH-property and $\mathcal{H}_{2}$ has $(\mathcal{G},v,\eta)$-strong EH-property
for some vectors $v,\eta$. 

From Theorem~\ref{product_theorem} we immediately obtain the following result that makes regular tournaments interesting in the context of the Erd\H{o}s-Hajnal Conjecture.

\begin{theorem}
\label{regular_theorem}
If $H$ is regular then it satisfies the Erd\"{o}s-Hajnal Conjecture.
\end{theorem}

In the next section we will give examples of new tournaments satisfying the Conjecture. All of them are regular.

We end this section with the proof of Theorem~\ref{product_theorem}.

\Proof
Take a $(v,\eta,c_{1},\lambda)$-$m$-sequence $\chi$, where $c_{1}>0$ is an arbitrary positive constant and $\lambda \leq \lambda_{0}(c_{1})$ ($\lambda_{0}(c_{1})$ as in the definition of the strong EH-property for $\mathcal{H}$). Since $\mathcal{H}_{1}$ has $(\mathcal{F},v,\eta)$-strong EH-property, by Theorem~\ref{strong_superstrong}, it also has $(\mathcal{F},v,\eta)$-super-strong EH-property. We can assume, without loss of generality, that some $H^{1}_{k} \in \mathcal{H}_{1}$ is $(f_{k},\frac{1}{2m})$-well-embedded in $\chi$, where:
$m = \sum_{i} \eta(i) + \sum_{j}[v(j)=0]$. Otherwise we get a $c_{2}$-strong pair for some constant $c_{2}>0$ and we are done. Let $V(H^{1}_{k})=\{h^{1}_{1},...,h^{1}_{|H_{k}|}\}$.
Let $l(\chi)=(T_{1},...,T_{m})$. Let $v_{1},...,v_{|H^{1}_{k}|}$ be the vertices of the embedding, where: $v_{j} \in T_{f_{k}(h^{1}_{j})}$ for $j=1,2,...,|H^{1}_{k}|$.
For each $T_{i}$ s.t. $v_{1},...,v_{|H^{1}_{k}|} \notin T_{i}$ denote: $R_{i} = \bigcap_{j=1,2,...,|H^{1}_{k}|} R^{j}_{i}$, where: $R^{j}_{i}$ is defined to be the set of the vertices of $T_{i}$ adjacent from $v_{j}$ if $f_{k}(h^{1}_{j}) < i$ and is defined to be the set of the vertices of $T_{i}$ adjacent to $v_{j}$ if  $f_{k}(h^{1}_{j}) > i$. From the fact that we have an $(f_{k},\frac{1}{2m})$-well-embedding, we conclude that $|R^{j}_{i}| \geq (1-\frac{1}{2m})|T_{i}|$, so: $|R_{i}| \geq \frac{1}{2}|T_{i}|$. 
Now we extract from $\chi$ a subsequence consisting of the sets $T_{i}$ s.t. $v_{1},...,v_{|H^{1}_{k}|} \notin T_{i}$. Then in that subsequence we replace each $T_{i}$ by $R_{i}$ to get a new $m$-sequence $\chi^{'}$. Note that $\mathcal{H}_{2}$ has $(\mathcal{G},v,\eta)$-strong EH-property. Thus we can assume without loss of generality that we have an embedding of some $H^{2}_{r} \in \mathcal{H}_{2}$ in $\chi^{'}$.
This embedding is defined by the vertices $w_{1},...,w_{|H^{2}_{r}|}$ s.t.  $w_{j} \in T_{g_{r}(w_{j})}$ for $j=1,2,...,|H^{2}_{r}|$ and some $g_{r} \in \mathcal{G}$.
Otherwise we get a $c_{2}$-strong pair for some constant $c_{2}>0$ and we are done. Combining this embedding with the embedding of $H^{1}_{k}$ we obtain an embedding of the tournament $(H^{1}_{k})^{f_{k}} \oplus (H^{2}_{r})^{g_{r}}$ from $\mathcal{H}_{1}^{\mathcal{F}} \oplus \mathcal{H}_{2}^{\mathcal{G}}$, 
defined by the set of vertices $\{v_{1},...,v_{|H^{1}_{k}|},w_{1},...,w_{|H^{2}_{r}|}\}$, where: $v_{j} \in T_{f_{k}(h^{1}_{j})}$ for $j=1,2,...,|H^{1}_{k}|$ and $w_{j} \in T_{g_{r}(w_{j})}$ for $j=1,2,...,|H^{2}_{r}|$.
That completes the proof. 
\bbox

\section{Applications}

In this section we give several applications of the strong EH-property and theorems proven by us in the previous sections. In particular we prove the Conjecture for all $6$-vertex tournaments but one and define new infinite families of tournaments satisfying the Conjecture. All those examples do not form the entire list of possible applications. One can define many other tournaments for which the Conjecture was open before but may be proven using our techniques.
It is worth to mention that the Conjecture for all prime tournaments for which it was known before follows from the techniques we present in this paper. 

We will give now few more technical definitions and observations that we will use in the Application section.

%Consider a $(v,\eta,c_{1},\lambda)$-$m$-sequence $\chi=(S_{1},...,S_{m})$ with %$l(\chi)=(T_{1},...,T_{k})$.
%Let $C>0$ be some constant.
%We say that a vertex $v \in S_{i}$ is \textit{$C$-smooth} if:
%\begin{itemize}
%\item $d(\{v\},S_{j}) \geq 1 - C\lambda$ for $j>i$ and
%\item $d(S_{j},\{v\}) \geq 1-  C\lambda$ for $j<i$.
%\end{itemize}
%We say that a set $\mathcal{S}$ is \textit{$C$-smooth} if every $s \in %\mathcal{S}$ is $C$-smooth. 

Let $D$ be a directed graph and let $V(D)=\{d_{1},...,d_{|D|}\}$. Let $\phi: V(D) \rightarrow \{1,2,...,k\}$ be an injective function. We say that $D$ is \textit{$(v,\eta,\phi)$-proper} if for an arbitrary $c>0$, $0 \leq \lambda \leq 1$ and $(v,\eta,c,\lambda)$-$m$-sequence $\chi$ with $l(\chi)=(T_{1},...,T_{k})$ there exists $c_{1}(v,\eta,c)$ such that: 

\begin{itemize}
\item $\exists v_{d_{1}},...,v_{d_{|H|}}$ such that $v_{d_{i}} \in T_{\phi(d_{i})}$ and: if $(d_{i},d_{j})$ is an edge then $(v_{d_{i}},v_{d_{j}})$ is an edge or
\item $\chi$ contains a $c_{1}$-strong pair.
\end{itemize}   

We have the following theorem:

\begin{theorem}
\label{directed_theorem}
Let $D_{1},...,D_{r}$ be directed graphs and assume that $D_{i}$ is $(v,\eta,\phi_{i})$-proper, where: $\phi_{i}(V(D_{i})) \cap \phi_{j}(V(D_{j})) = \emptyset$ for $i \neq j$. 
Let $V(D_{i})=\{d^{i}_{1},...,d^{i}_{|D_{i}|}\}$.
Let $\chi$ be a $(v,\eta,c,\lambda)$-$m$-sequence with $l(\chi)=(T_{1},...,T_{k})$, where $c>0$ is an arbitrary constant and $\lambda > 0$ is small enough. 
Then there exists $c_{1}>0$ s.t.:
\begin{itemize}
\item there exist vertices: $v_{d^{1}_{1}},...,v_{d^{r}_{|D_{r}|}}$ s.t.:
\begin{itemize}
\item $v_{d^{i}_{j}} \in T_{\phi_{i}(j)}$, 
\item if $(d^{k}_{i},d^{k}_{j})$ is an edge then $(v_{d^{k}_{i}},v_{,d^{k}_{j}})$ is an edge, 
\item there are no edges $(v_{d^{k_{1}}_{j}},v_{d^{k_{2}}_{i}})$ if $k_{1} \neq k_{2}$, $\phi_{k_{1}}(j) > \phi_{k_{2}}(i)$, or
\end{itemize} 
\item $\chi$ contains a $c_{1}$-strong pair.
\end{itemize}
\end{theorem}

\Proof
The proof is by induction on $r$. For $r=1$ the result is trivial. Thus assume that $r>1$. Let $C=2k$.
As in the proof of Theorem~\ref{strong_superstrong}, we can replace $\chi$ by a $(v,\eta,c^{'},\lambda^{'})$-strong $m$-sequence $\chi^{'}$, where: $c^{'}=c(1-\frac{k}{C})$, $\lambda^{'}=\frac{\frac{C\lambda}{c^{2}}}{1-\frac{k}{C}}$. 
From the fact that $D_{1}$ is $(\chi,\phi_{1},c_{1})$-proper we can conclude that we either have a strong pair in $\chi^{'}$ and we are done or there exists an embedding $v_{d^{1}_{1}},...,v_{d^{1}_{|D_{1}|}}$ as in the definition of the $(\chi,\phi_{1},c_{1})$-properness. Thus assume the latter. For every $i$ s.t. $v_{d^{1}_{1}},...,v_{d^{1}_{|D_{1}|}} \notin T_{i}$ we define $R_{i}=\bigcap_{j=1,...,|D_{1}|} R^{j}_{i}$, where $R^{j}_{i}$ is the subset of the vertices of $T_{i}$ that are adjacent from $v_{d^{1}_{j}}$ if $\phi_{1}(d^{1}_{j}) < i$ and are adjacent to $v_{d^{1}_{j}}$ if $\phi_{1}(d^{1}_{j}) > i$. Let $\chi^{''}$ be a $m$-sequence obtained from $\chi^{'}$ by replacing every $T_{i}$ s.t. $v_{d^{1}_{1}},...,v_{d^{1}_{|D_{1}|}} \notin T_{i}$ by $R_{i}$.
It is easy to see that by induction for $\lambda>0$ small enough the following is true:
\begin{itemize}
\item there exists a strong pair in $\chi^{''}$ or
\item there exists an embedding in $\chi^{''}$ as in the statement of the theorem, for a family: $\{D_{2},...,D_{r}\}$.
\end{itemize}

If the former is true we are done. Thus assume that the latter is true.
Then, according to the definition of the sets $R_{i}$, we are also done since we can combine the embedding $\{v^{1}_{1},...,v^{1}_{|D_{1}|}\}$ with the embedding we have just mentioned and get an embedding for $\{D_{1},...,D_{r}\}$, as in the statement of the theorem. 
\bbox

\subsection{Galaxies}

We start by proving that all galaxies satisfy the Conjecture. The Conjecture was first proven for galaxies in \cite{bcc}. However we would like to show that it easily follows from the theorems concerning strong EH-property.

\begin{theorem}
Every galaxy satisfies the Erd\"{o}s-Hajnal Conjecture.
\end{theorem}

\Proof

Take a left star $H$ of $l$ leaves with $V(H)=\{h_{1},...,h_{l+1}\}$. Assume that $h_{1}$ is a center. 
Let $v=(0,1)$ and $\eta = (l)$. Let $f:V(H) \rightarrow \mathbb{N}$ be a function defined as follows: $f(h_{j})=j$.
%???
It can be shown that $H$ has $(\{f\},v,\eta)$-strong EH-property (this is what in fact was proven in \cite{bcc}, though in that paper the proof wasnt interpreted as an application of the much more general method).
The similar result may be proven for the right star.
Let $\mathcal{H}$ be a family of galaxies. Note that $\mathcal{H}=\hat{\mathcal{H}_{b}}$, where $\mathcal{H}_{b}$ is the set of left and right stars (i.e. the family of galaxies may be obtained from the set of left and right stars by taking EH-extensions and applying operation $\oplus$ on pairs of tournaments).
That, according to Theorem~\ref{extension_theorem}, completes the proof.
\bbox

Our methods may be also used to prove the Conjecture for a larger family of tournaments that contains galaxies, called \textit{constellations} (see: \cite{choromanski2}). However we will skip the proof here. 

All the remaining examples in this section will be families of tournaments that satisfy the Conjecture but for which the Conjecture was open before.

\subsection{Right pseudogalaxies}

Let us define the family of tournaments called \textit{right pseudogalaxies}.
Let $H = H_{l}^{\rho_{1}} \oplus G^{\rho_{2}}$, where $H_{l},G,\rho_{1},\rho_{2}$ are defined as follows. $V(H_{l})=\{h_{1},...,h_{7}\}$. Under ordering $(h_{1},...,h_{7})$ the backward edges of $H$ are: $\{(h_{4},h_{1}),(h_{6},h_{2}),(h_{6},h_{1}),(h_{5},h_{3}),(h_{7},h_{5})\}$.
Let $G$ be a galaxy with $V(G)=\{g_{1},...,g_{r}\}$, where: $(g_{1},...,g_{r})$ is a galaxy ordering. Let $1 \leq s \leq r$ be some constant and assume that $G$ does not have a star with one leaf in $\{g_{1},...,g_{s}\}$ and the other in $\{g_{s+1},...,g_{r}\}$. Let $\rho_{1}:V(H_{l}) \rightarrow \mathbb{N}$ be a function defined as: $\rho_{1}(h_{i})=i$ for $i=1,...,6$, $\rho_{1}(h_{7})=6+s+1$.
Let $\rho_{2}:V(G) \rightarrow \mathbb{N}$ be a function defined as: $\rho_{2}(g_{u})=6+u$ for $u=1,...,s$ and $\rho_{2}(g_{u})=7+u$ for $u=s+1,...,r$.
then we say that $H$ is a \textit{right pseudogalaxy}. The name comes from the fact that $H$ is created by combining a galaxy with another tournament using $\oplus$ operator and the galaxy is "attached" from the right.

\begin{theorem}
\label{right_pseudogalaxy_theorem}
Every right pseudogalaxy $H$ satisfies the Erd\H{o}s-Hajnal Conjecture.
\end{theorem}

\Proof
Let $H=H_{l}^{\rho_{1}} \oplus G^{\rho_{2}}$.
Consider the following vectors: $v=(0,0,1,0,0,0,0,0,0,0,0), \eta=(3)$. 
Consider 4 tournaments: $H_{1},H_{2},H_{3},H_{4}$ with:
\begin{itemize}
\item $V(H_{1})=\{v^{1}_{1},v^{1}_{3},v^{1}_{5},v^{1}_{9},v^{1}_{4},v^{1}_{7},v^{1}_{12}\}$,
\item $V(H_{2})=\{v^{2}_{1},v^{2}_{3},v^{2}_{9},v^{2}_{2},v^{2}_{8},v^{2}_{7},v^{2}_{12}\}$,
\item $V(H_{3})=\{v^{3}_{1},v^{3}_{5},v^{3}_{9},v^{3}_{2},v^{3}_{8},v^{3}_{7},v^{3}_{12}\}$,
\item $V(H_{4})=\{v^{4}_{4},v^{4}_{7},v^{4}_{12},v^{4}_{6},v^{4}_{11},v^{4}_{10},v^{4}_{13}\}$.
\end{itemize}

Let the sets of backward edges $B(H_{1}),B(H_{2}),B(H_{3}),B(H_{4})$ of $H_{1},H_{2},H_{3},H_{4}$ under orderings: $(v^{1}_{1},v^{1}_{3},v^{1}_{5},v^{1}_{9},v^{1}_{4},v^{1}_{7},v^{1}_{12}),
(v^{2}_{1},v^{2}_{3},v^{2}_{9},v^{2}_{2},v^{2}_{8},v^{2}_{7},v^{2}_{12}),
(v^{3}_{1},v^{3}_{5},v^{3}_{9},v^{3}_{2},v^{3}_{8},v^{3}_{7},v^{3}_{12}),
(v^{4}_{4},v^{4}_{7},v^{4}_{12},v^{4}_{6},v^{4}_{11},v^{4}_{10},v^{4}_{13})$ be respectively:

\begin{itemize}
\item $B(H_{1})=\{(v^{1}_{5},v^{1}_{1}),(v^{1}_{9},v^{1}_{3}),(v^{1}_{9},v^{1}_{1}),(v^{1}_{7},v^{1}_{4}),(v^{1}_{12},v^{1}_{7})\}$, 
\item $B(H_{2})=\{(v^{2}_{9},v^{2}_{3}),(v^{2}_{9},v^{2}_{1}),(v^{2}_{3},v^{2}_{1}),(v^{2}_{8},v^{2}_{2}),(v^{2}_{12},v^{2}_{7})\}$,
\item $B(H_{3})=\{(v^{3}_{9},v^{3}_{5}),(v^{3}_{9},v^{3}_{1}),(v^{3}_{5},v^{3}_{1}),(v^{3}_{8},v^{3}_{2}),(v^{3}_{12},v^{3}_{7})\}$,
\item $B(H_{4})=\{(v^{4}_{12},v^{4}_{7}),(v^{4}_{12},v^{4}_{4}),(v^{4}_{7},v^{4}_{4}),(v^{4}_{11},v^{4}_{6}),(v^{4}_{13},v^{4}_{10})\}$.
\end{itemize}

Consider the family $\mathcal{F}$ of functions defined as follows: $\mathcal{F} = \{f_{1},f_{2},f_{3},f_{4}\}$, where: $f_{1}:V(H_{1}) \rightarrow \{1,2,...,13\}$,  $f_{2}:V(H_{2}) \rightarrow \{1,2,...,13\}$,  $f_{3}:V(H_{3}) \rightarrow \{1,2,...,13\}$,  $f_{4}:V(H_{4}) \rightarrow \{1,2,...,13\}$.
Let $f_{i}(v^{i}_{j})=j$. 
The crucial observation now is that $\mathcal{H}=\{H_{1},H_{2},H_{3},H_{4}\}$ has $(\mathcal{F},v,\eta)$-strong EH-property. 

\begin{theorem}
\label{sEH_right}
The family $\mathcal{H}=\{H_{1},H_{2},H_{3},H_{4}\}$ has $(\mathcal{F},v,\eta)$-strong EH-property. 
\end{theorem}

\Proof
Let $D_{1}$ be a directed graph with $V(D_{1})=\{d^{1}_{1},d^{1}_{3},d^{1}_{5},d^{1}_{9}\}$ and $E(D_{1})=\{(d^{1}_{5},d^{1}_{1}),(d^{1}_{9},d^{1}_{3}),(d^{1}_{9},d^{1}_{1})\}$.
Let $\phi_{1}: V(D_{1}) \rightarrow \{1,...,13\}$ be defined as follows:
$\phi_{1}(d^{1}_{j}) = j$. Note that $D_{1}$ is $(v,\eta,\phi_{1})$-proper.
Indeed, take a $(v,\eta,c,\lambda)$-$m$-sequence $\tau$ with $l(\tau)=(V_{1},...,V_{13})$. Let:
\begin{itemize}
\item  $V_{1}^{'}=\{v \in V_{1}: \exists x_{v} \in V_{5}$ s.t. $(x_{v},v)\}$,
\item  $V_{9}^{'}=\{v \in V_{9}: \exists x_{v} \in V_{3}$ s.t. $(v,x_{v})\}$.
\end{itemize}
We may assume that $|V_{1}^{'}| \geq \frac{|V_{1}|}{2}$. Otherwise we get a $\frac{c}{2}$-strong pair $(V_{1} \backslash V_{1}^{'},V_{5})$ (note that $V_{1} \backslash V_{1}^{'}$ is complete to $V_{5}$). Similarly, we may assume that 
$|V_{9}^{'}| \geq \frac{|V_{9}|}{2}$. Now notice that we may assume that there exist vertices $x,y$ s.t. $x \in V_{1}^{'}$, $y \in V_{9}^{'}$ and $y$ is adjacent to $x$. Otherwise we have: $(V_{1}^{'},V_{9}^{'})$ is a $\frac{c}{2}$-strong pair. Let $x^{'}$ be an arbitrary vertex in $V_{5}$ adjacent to $x$ and let $y^{'}$ be an arbitrary vertex in $V_{3}$ adjacent from $y$. Then $\{x,y^{'},x^{'},y\}$ is the embedding from the definition of the $(v,\eta,\phi_{1})$-properness of $D_{1}$. That completes the proof that $D_{1}$ is $(v,\eta,\phi_{1})$-proper. 
Let $D_{2}$ be a directed graph with $V(D_{2})=\{d^{2}_{4},d^{2}_{7},d^{2}_{12}\}$ and $E(D_{2})=\{(d^{2}_{7},d^{2}_{4}),(d^{2}_{12},d^{2}_{7})\}$. Let $\phi_{2}: V(D_{2}) \rightarrow \{1,...,13\}$ be defined as follows:
$\phi_{2}(d^{2}_{j}) = j$. By the completely analogous analysis we can conclude that $D_{2}$ is $(v,\eta,\phi_{2})$-proper.
Let $D_{3}$ be a directed graph with $V(D_{3})=\{d^{3}_{2},d^{3}_{8}\}$ and $E(D_{3})=\{(d^{3}_{8},d^{3}_{2})\}$. Let $\phi_{3}: V(D_{3}) \rightarrow \{1,...,13\}$ be defined as follows:
$\phi_{3}(d^{3}_{j}) = j$. 
Let $D_{4}$ be a directed graph with $V(D_{4})=\{d^{4}_{6},d^{4}_{11}\}$ and $E(D_{4})=\{(d^{4}_{11},d^{4}_{6})\}$. Let $\phi_{4}: V(D_{4}) \rightarrow \{1,...,13\}$ be defined as follows:
$\phi_{4}(d^{4}_{j}) = j$. 
Let $D_{5}$ be a directed graph with $V(D_{5})=\{d^{5}_{13},d^{5}_{10}\}$ and $E(D_{5})=\{(d^{5}_{13},d^{5}_{10})\}$. Let $\phi_{5}: V(D_{5}) \rightarrow \{1,...,13\}$ be defined as follows:
$\phi_{5}(d^{5}_{j}) = j$. 
Trivially, $D_{i}$ is $(v,\eta,\phi_{i})$-proper for $i=3,4,5$.
Now we can apply Theorem~\ref{directed_theorem} for directed graphs: $D_{1},...,D_{5}$ and assume, without loss of generality, that as an outcome we get an embedding from the statement of the theorem. It remains to notice that this embedding gives an embedding from the definition of the $(\mathcal{F},v,\eta)$-strong EH-property of $\mathcal{H}$.

\bbox

Since the family $\mathcal{H}$  has $(\mathcal{F},v,\eta)$-strong EH-property, for every EH-extension $(\mathcal{F}^{1},v^{1},\eta^{1})$ of $(\mathcal{F},v,\eta)$ the family $\mathcal{H}$ has $(\mathcal{F}^{1},v^{1},\eta^{1})$-strong EH-property.
If $H_{g}$ is an arbitrary galaxy then from the previous example we know that $\{H_{g}\}$ has $(\mathcal{F}^{2},v^{2},\eta^{2})$-strong EH-property for appropriate vectors $v^{2},\eta^{2}$ and one-element family $\mathcal{F}^{2}=\{f\}$, where $f$ is an appropriate function inducing a galaxy ordering of the vertices of $H_{g}$. Now, from Theorem~\ref{product_theorem} we know that if $<\mathcal{F}^{1},\mathcal{F}^{2}>=0$, then $\mathcal{H}^{\mathcal{F}^{1}} \oplus \{H_{g}\}^{\mathcal{F}^{2}}$ has strong EH-property.
It remains to notice that there exists an extension $(\mathcal{F}^{1},v^{1},\eta^{1})$, a triple $(\mathcal{F}^{2},v^{2},\eta^{2})$  and a galaxy $H_{g}$ such that we have:\\ 

every tournament in $\mathcal{H}^{\mathcal{F}^{1}} \oplus \{H_{g}\}^{\mathcal{F}^{2}}$ contains $H_{l}^{\rho_{1}} \oplus G^{\rho_{2}}$ or $H_{t}^{\rho_{3}} \oplus G^{\rho_{2}}$ as a subtournament, \\

where:
$V(H_{t})=\{k_{1},...,k_{7}\}$, $\rho_{3}$ is defined as follows:

\begin{itemize}
\item $\rho_{3}(k_{1})=3$,
\item $\rho_{3}(k_{2})=1$,
\item $\rho_{3}(k_{3})=5$,
\item $\rho_{3}(k_{4})=2$,
\item $\rho_{3}(k_{5})=4$,
\item $\rho_{3}(k_{6})=6$,
\item $\rho_{3}(k_{7})=7+s$,
\end{itemize} 

and the set of backward edges of $H_{t}$ under an ordering induced by $\rho_{3}$ is: $$\{(k_{1},k_{2}),(k_{6},k_{1}),(k_{6},k_{2}),(k_{3},k_{4}),(k_{7},k_{5})\}.$$

We complete the proof noticing that  $H_{l}^{\rho_{1}} \oplus G^{\rho_{2}}=H_{t}^{\rho_{3}} \oplus G^{\rho_{2}}$. Thus every right pseudogalaxy is regular. 
\bbox

\subsection{$C_{5}$-chains and the Conjecture for tournaments on $6$ vertices}

Consider the following family $\mathcal{H}$ of tournaments.

Every element $H$ of $\mathcal{H}$ is of the form: $H_{l}^{\rho_{1}} \oplus P_{2t+1}^{\rho_{2}}$, where:
\begin{itemize}
\item $V(H_{l})=\{h_{1},...,h_{6}\}$,
\item $\rho_{1}(h_{i}) = i$ for $i=1,2,...,5$,
\item $\rho_{1}(h_{6}) = 7$,
\item the set of backward edges of $H_{l}$ under an ordering $(h_{1},...,h_{6})$ is $\{(h_{4},h_{1}),(h_{5},h_{2}),(h_{5},h_{1}),(h_{6},h_{3})\}$,
\item $V(P_{2t+1})=\{p_{1},...,p_{2t+1}\}$,
\item $\rho_{2}(p_{2i+1})=2i+6$, $\rho_{2}(p_{2i})=2i+7$,
\item the set of backward edges of $P_{2t+1}$ under an ordering $(p_{1},...,p_{2t+1})$ is $\{(p_{2i},p_{2i-1}): i=1,2,...\}$.
\end{itemize}

We call $H$ a \textit{$C_{5}$-chain}. The $C_{5}$-part of the name comes from the fact that every $C_{5}$-chain contains $C_{5}$ as a subtournament.

\begin{theorem}
\label{chainresult}
Every $C_{5}$-chain $H$ satisfies the Erd\H{o}s-Hajnal Conjecture. 
\end{theorem}

\Proof
Let $v=(0,0,1,0,...,0)$ be a $(2t+7)$-element vector.
Let $\eta=(3)$. Let us consider three tournaments: $H_{1},H_{2},H_{3}$ such that:
\begin{itemize}
\item $V(H_{1})= \{v^{1}_{1},v^{1}_{3},v^{1}_{5},v^{1}_{8},v^{1}_{4},v^{1}_{10},v^{1}_{13},v^{1}_{9},v^{1}_{2t+6}\} \cup \{v^{1}_{4i+3},v^{1}_{4i}:i \geq 3, 4i \leq 2t+4\} \cup  \{v^{1}_{4i+1},v^{1}_{4i-2}:i \geq 4, 4i \leq 2t+6\}$,
\item $V(H_{2})=\{v^{2}_{1},v^{2}_{3},v^{2}_{8},v^{2}_{2},v^{2}_{6},v^{2}_{11},v^{2}_{7},v^{2}_{13},v^{2}_{9}\}\cup \{v^{2}_{5i+1},v^{2}_{5i-4}: i \geq 3, 5i \leq 2t+6\} \cup \{v^{2}_{5i+3},v^{2}_{5i-2}: i \geq 3, 5i \leq 2t+4\}$,
\item 
$V(H_{3}) = \{v^{3}_{1},v^{3}_{5},v^{3}_{8},v^{3}_{2},v^{3}_{6},v^{3}_{11},v^{3}_{7},v^{3}_{13},v^{3}_{9}\} \cup \{v^{3}_{5i+1},v^{3}_{5i-4}: i \geq 3, 5i \leq 2t+6\} \cup \{v^{3}_{5i+3},v^{3}_{5i-2}: i \geq 3, 5i \leq 2t+4\}$,
\item the set of backward edges $B(H_{1})$ of $H_{1}$ under an ordering of the vertices $v^{1}_{j}$ of $H_{1}$ induced by an increasing value of $j$ is:
$B(H_{1}) = \{(v^{1}_{5},v^{1}_{1}),(v^{1}_{8},v^{1}_{3}),(v^{1}_{8},v^{1}_{1}),(v^{1}_{13},v^{1}_{9}),(v^{1}_{15},v^{1}_{12})\} \cup \{(v^{1}_{4i+3},v^{1}_{4i}):i \geq 3, 4i \leq 2t+4\} \cup  \{(v^{1}_{4i+1},v^{1}_{4i-2}):i \geq 4, 4i \leq 2t+6\}$,

\item the set of backward edges $B(H_{2})$ of $H_{2}$ under an ordering of the vertices $v^{2}_{j}$ of $H_{2}$ induced by an increasing value of $j$ is:
$B(H_{2}) = \{(v^{2}_{3},v^{2}_{1}),(v^{2}_{8},v^{2}_{1}),(v^{2}_{8},v^{2}_{3}),(v^{2}_{6},v^{2}_{2}),(v^{2}_{11},v^{2}_{7}),(v^{2}_{13},v^{2}_{9})\} \\ \cup \{(v^{2}_{5i+1},v^{2}_{5i-4}): i \geq 3, 5i \leq 2t+6\} \cup \{(v^{2}_{5i+3},v^{2}_{5i-2}): i \geq 3, 5i \leq 2t+4\}$,

\item the set of backward edges $B(H_{3})$ of $H_{3}$ under an ordering of the vertices $v^{3}_{j}$ of $H_{3}$ induced by an increasing value of $j$ is:
$B(H_{3}) = \{(v^{3}_{5},v^{3}_{1}),(v^{3}_{8},v^{3}_{1}),(v^{3}_{8},v^{3}_{5}),(v^{3}_{6},v^{3}_{2}),(v^{3}_{11},v^{3}_{7}),(v^{3}_{13},v^{3}_{9})\}
\\ \cup \{(v^{3}_{5i+1},v^{3}_{5i-4}): i \geq 3, 5i \leq 2t+6\} \cup \{(v^{3}_{5i+3},v^{3}_{5i-2}): i \geq 3, 5i \leq 2t+4\}.
$
\end{itemize}

Let us consider the following family $\mathcal{F}$ of functions: $\mathcal{F}=\{f_{1},f_{2},f_{3}\}$, where: 
$f_{1}:V(H_{1}) \rightarrow \{1,...,2t+7\}$, 
$f_{2}:V(H_{2}) \rightarrow \{1,...,2t+7\}$, 
$f_{3}:V(H_{3}) \rightarrow \{1,...,2t+7\}$ 
and $f_{i}(v^{i}_{j}) = j$.

%???
As in the previous section, we can note that $\mathcal{H}=\{H_{1},H_{2},H_{3}\}$ has $(\mathcal{F},v,\eta)$-strong EH-property. 
The proof is completely analogous to the related one from the previous section and the set of directed graphs used there is a subset of the set of directed graphs used in the corresponding proof for right pseudogalaxies.

Note that $H_{1}=H_{2}=H_{3}=H$. Thus $H$ is regular and the proof is completed.
\bbox

Let $H_{D}$ be a tournament with $V(H_{D})=\{1,2,3,4,5,6\}$, where under ordering $(1,2,3,4,5,6)$ of its vertices the set of backward edges is: ${(4,1),(6,3),(5,2),(6,1)}$.

Now we will prove the corollary of Theorem~\ref{chainresult}, namely:
\begin{theorem}
If $H$ is a tournament on six vertices not isomorphic to $H_{D}$ then it satisfies the Erd\H{o}s-Hajnal Conjecture.
\end{theorem}

\Proof
It was proven in \cite{bcc} that every tournament on at most $5$ vertices satisfies the Erd\H{o}s-Hajnal Conjecture.
Let us define two tournaments: $H_{a},H_{b}$ as follows: $V(H_{a})=\{1,...,6\}$, $V(H_{b})=\{1,...,6\}$. The set $B(H_{a})$ of backward edges of $H_{a}$ under ordering $(1,...6)$ of its vertices is: $\{(4,1),(5,2),(5,1),\\(6,3)\}$. The set $B(H_{b})$ of backward edges of $H_{b}$ under ordering $(1,2,3,4,6,5)$ of its vertices is the same as for $H_{a}$. It was proven (see: \cite{mc}) that every tournament on six vertices is either a galaxy, may be obtained from smaller tournaments by a substitution procedure or is isomorphic to $H_{a}$, $H_{b}$ or $H_{D}$.
Now it remains to note that $H_{b}$ is a subtournament of some right pseudogalaxy and $H_{a}$ is a subtournament of some $C_{5}$-chain. Thus the result follows from theorems:~\ref{right_pseudogalaxy_theorem} and ~\ref{chainresult}.
\bbox

\subsection{Left pseudogalaxies and $C_{5}$-co-chains}

Tournament $H$ is a \textit{left pseudogalaxy} if it is a complement of some right pseudogalaxy. Similarly, tournament $H$ is a $C_{5}$-co-chain if it is a complement of some $C_{5}$ chain. Note that if a tournament $H_{2}$ is the complement of the tournament $H_{1}$ then the graph of backward edges of $H_{2}$ is the mirror image of the graph of backward edges of $H_{1}$.
Thus the proofs of theorems:~\ref{right_pseudogalaxy_theorem} and~\ref{chainresult} may be almost exactly repeated and we get the following theorem:

\begin{theorem}
\label{mirror_image_theorem}
If $H$ is a left pseudogalaxy or a $C_{5}$-co-chain then it satisfies the Erd\"{o}s-Hajnal Conjecture. 
\end{theorem}

\subsection{$\alpha$-galaxies}

Tournament $H$ is a $\alpha$-galaxy if $H=H_{s}^{\rho_{1}} \oplus G^{\rho_{2}}$, where: 

\begin{itemize}
\item $V(H_{s}) = \{h_{1},,,.,h_{7}\}$,
\item the set $B(H_{s})$ of backward edges of $H_{S}$ under ordering $(h_{1},...,h_{7})$ is of the form:\\ $\{(h_{4},h_{1}),(h_{6},h_{4}),(h_{6},h_{3}),(h_{7},h_{5})\}$,
\item $G$ is a galaxy with $V(G)=\{g_{1},...,g_{r}\}$ and galaxy ordering of its vertices: $(g_{1},...,g_{r})$,
\item  there exists: $0 \leq r_{1} \leq r$ s.t.
for $\mathcal{R}_{1}=\{g_{1},...,g_{r_{1}}\}$, $\mathcal{R}_{2} = \{g_{r_{1}+1},...,g_{r}\}$ the following is true: no star of $G$ has one leaf in $\mathcal{R}_{1}$ and the other one in $\mathcal{R}_{2}$, 
\item for $r_{1}$ given above there exists $r_{2},r_{3} \geq 0$ s.t. $r=r_{1}+r_{2}+r_{3}$ and: $\rho_{2}(g_{i}) = 5+i$ for $i=1,...,r_{1}$, $\rho_{2}(g_{i})=6+i$ for $i=r_{1}+1,...,r_{1}+r_{2}$, $\rho_{2}(g_{i})=7+i$ for $i=r_{1}+r_{2}+1,...,r$,
\item $\rho_{1}(h_{i})=i$ for $i=1,...,5$, $\rho_{1}(h_{6})=6+r_{1}$, $\rho_{1}(h_{7})=7+r_{1}+r_{2}$.
 
\end{itemize}

\begin{theorem}
\label{alfa_theorem}
If $H$ is a $\alpha$-galaxy then it satisfies the Erd\H{o}s-Hajnal Conjecture.
\end{theorem}

\Proof
Let $v=(0,1,0,0,0,0,0,1)$, $\eta=(3,1)$. Let us consider three tournaments: $H_{1},H_{2},H_{3}$ s.t.:

\begin{itemize}
\item $V(H_{1}) = \{v^{1}_{2},v^{1}_{3},v^{1}_{7},v^{1}_{4},v^{1}_{9},v^{1}_{8},v^{1}_{10}\}$,
\item $V(H_{2}) = \{v^{2}_{1},v^{2}_{6},v^{2}_{2},v^{2}_{7},v^{2}_{9},v^{2}_{8},v^{2}_{10}\}$,
\item $V(H_{3}) = \{v^{3}_{3},v^{3}_{5},v^{3}_{4},v^{3}_{7},v^{3}_{9},v^{3}_{8},v^{3}_{10}\}$,
\item the set $B(V(H_{1}))$ of backward edges of $H_{1}$ under an ordering of the vertices $v^{1}_{j}$ of $H_{1}$ induced by an increasing value of $j$ is: 
$B(V(H_{1})) = \{(v^{1}_{9},v^{1}_{7}),(v^{1}_{9},v^{1}_{4}),(v^{1}_{7},v^{1}_{2}),(v^{1}_{10},v^{1}_{8})\}$,
\item the set $B(V(H_{2}))$ of backward edges of $H_{2}$ under an ordering of the vertices $v^{2}_{j}$ of $H_{1}$ induced by an increasing value of $j$ is: 
$B(V(H_{2})) = \{(v^{2}_{6},v^{2}_{1}),(v^{2}_{9},v^{2}_{7}),(v^{2}_{9},v^{2}_{2}),(v^{2}_{7},v^{2}_{2}),(v^{2}_{10},v^{2}_{8})\}$,
\item the set $B(V(H_{3}))$ of backward edges of $H_{2}$ under an ordering of the vertices $v^{3}_{j}$ of $H_{1}$ induced by an increasing value of $j$ is: 
$B(V(H_{3})) = \{(v^{3}_{5},v^{3}_{3}),(v^{3}_{9},v^{3}_{4}),(v^{3}_{9},v^{3}_{7}),(v^{3}_{7},v^{3}_{4}),(v^{3}_{10},v^{3}_{8})\}$.
\end{itemize}

Let $\mathcal{F}=\{f_{1},f_{2},f_{3}\}$ be the family of functions, where: $f_{i}(v^{i}_{j})=j$.
Let $\mathcal{H}=\{H_{1},H_{2},H_{3}\}$.
%??? tiger
Now we observe that $\mathcal{H}=\{H_{1},H_{2},H_{3}\}$ has $(\mathcal{F},v,\eta)$-strong EH-property. 

\begin{theorem}
The family $\mathcal{H}=\{H_{1},H_{2},H_{3}\}$ has $(\mathcal{F},v,\eta)$-strong EH-property.
\end{theorem}

\Proof
Let $D_{1}$ be a directed graph with $V(D_{1})=\{d^{1}_{2},d^{1}_{4},d^{1}_{7},d^{1}_{9}\}$ and $E(D_{1})=\{(d^{1}_{7},d^{1}_{2}),(d^{1}_{9},d^{1}_{7}),(d^{1}_{9},d^{1}_{4})\}$.
Let $\phi_{1}: V(D_{1}) \rightarrow \{1,...,10\}$ be defined as follows:
$\phi_{1}(d^{1}_{j}) = j$. Note that $D_{1}$ is $(v,\eta,\phi_{1})$-proper.
Indeed, take a $(v,\eta,c,\lambda)$-$m$-sequence $\tau$ with $l(\tau)=(V_{1},...,V_{10})$. 
Let:
\begin{itemize}
\item  $V_{7}^{'}=\{v \in V_{7}: \exists x_{v} \in V_{2}$ s.t. $(v,x_{v})\}$.
\end{itemize}
Note that we can assume that $|V_{7}^{'}| \geq \frac{|V_{7}|}{2}$. Otherwise we
get a $\frac{c}{2}$-strong pair $(V_{7} \backslash V_{7}^{'},V_{2})$.
Let:
\begin{itemize}
\item  $V_{9}^{'}=\{v \in V_{9}: \exists x_{v} \in V_{7}^{'}$ s.t. $(v,x_{v})\}$.
\end{itemize}
Note that we can assume that $|V_{9}^{'}| \geq \frac{|V_{9}|}{2}$. Otherwise, since $|V_{7}^{'}| \geq \frac{|V_{7}|}{2}$, we
get a $\frac{c}{2}$-strong pair $(V_{7}^{'},V_{9} \backslash V_{9}^{'})$.
Let:
\begin{itemize}
\item  $V_{9}^{''}=\{v \in V_{9}^{'}: \exists x_{v} \in V_{4}$ s.t. $(v,x_{v})\}$.
\end{itemize}
Note that we can assume that $|V_{9}^{''}| \geq \frac{|V_{9}^{'}|}{2}$. Otherwise we get a $\frac{c}{4}$-strong pair $(V_{9}^{''} \backslash V_{9}^{'},V_{4})$.
Take a vertex $x \in V_{9}^{''}$. From the definition of $V_{9}^{''}$ we know that there exist $a \in V_{2}, b \in V_{4}, c \in V_{7}$ s.t. $x$ is adjacent to $b$ and $c$ and $c$ is adjacent to $a$. But then $\{a,b,c,x\}$ is the embedding from the definition of the $(v,\eta,\phi_{1})$-properness of $D_{1}$ and that completes the proof of the fact that $D_{1}$ is  $(v,\eta,\phi_{1})$-proper.
Let $D_{2}$ be a directed graph with $V(D_{2})=\{d^{2}_{6},d^{2}_{1}\}$ and $E(D_{2})=\{(d^{2}_{6},d^{2}_{1})\}$. Let $\phi_{2}: V(D_{2}) \rightarrow \{1,...,10\}$ be defined as follows:
$\phi_{2}(d^{2}_{j}) = j$. 
Let $D_{3}$ be a directed graph with $V(D_{3})=\{d^{3}_{5},d^{3}_{3}\}$ and $E(D_{3})=\{(d^{3}_{5},d^{3}_{3})\}$. Let $\phi_{3}: V(D_{3}) \rightarrow \{1,...,10\}$ be defined as follows:
$\phi_{3}(d^{3}_{j}) = j$. 
Let $D_{4}$ be a directed graph with $V(D_{4})=\{d^{4}_{10},d^{4}_{8}\}$ and $E(D_{4})=\{(d^{4}_{10},d^{4}_{8})\}$. Let $\phi_{4}: V(D_{4}) \rightarrow \{1,...,10\}$ be defined as follows:
$\phi_{4}(d^{4}_{j}) = j$. 
Trivially, $D_{i}$ is $(v,\eta,\phi_{i})$-proper for $i=2,3,4$.
Now we can apply Theorem~\ref{directed_theorem} for directed graphs: $D_{1},...,D_{4}$ and assume, without loss of generality, that as an outcome we get an embedding from the statement of the theorem. It remains to notice that this embedding gives an embedding from the definition of the $(\mathcal{F},v,\eta)$-strong EH-property of $\mathcal{H}$. 

\bbox

From what we have just proven we know that for every EH-extension $(\mathcal{F}^{1},v^{1},\eta^{1})$ of $(\mathcal{F},v,\eta)$ the family $\mathcal{H}$ has $(\mathcal{F}^{1},v^{1},\eta^{1})$-strong EH-property.
If $H_{g}$ is an arbitrary galaxy then we know that $\{H_{g}\}$ has $(\mathcal{F}^{2},v^{2},\eta^{2})$-strong EH-property for appropriate vectors $v^{2},\eta^{2}$ and one-element family $\mathcal{F}^{2}=\{f\}$, where $f$ is an appropriate function inducing a galaxy ordering of the vertices of $H_{g}$. Now, from Theorem~\ref{product_theorem} we know that if $<\mathcal{F}^{1},\mathcal{F}^{2}>=0$, then $\mathcal{H}^{\mathcal{F}^{1}} \oplus \{H_{g}\}^{\mathcal{F}^{2}}$ has strong EH-property.
It remains to notice that there exists an extension $(\mathcal{F}^{1},v^{1},\eta^{1})$, a triple $(\mathcal{F}^{2},v^{2},\eta^{2})$ and a galaxy $H_{g}$ such that we have:\\ 

every tournament in $\mathcal{H}^{\mathcal{F}^{1}} \oplus \{H_{g}\}^{\mathcal{F}^{2}}$ contains $H_{s}^{\rho_{1}} \oplus G^{\rho_{2}}$ or $H_{t}^{\rho_{3}} \oplus G^{\rho_{2}}$ as a subtournament, \\

where:
$V(H_{t})=\{k_{1},...,k_{7}\}$, $\rho_{3}$ is defined as follows:

\begin{itemize}
\item $\rho_{3}(k_{1})=3$,
\item $\rho_{3}(k_{2})=1$,
\item $\rho_{3}(k_{3})=4$,
\item $\rho_{3}(k_{4})=2$,
\item $\rho_{3}(k_{5})=5$,
\item $\rho_{3}(k_{6})=6+r_{1}$,
\item $\rho_{3}(k_{7})=7+r_{1}+r_{2}$,
\end{itemize} 

and the set of backward edges of $H_{t}$ under an ordering induced by $\rho_{3}$ is: $$\{(k_{1},k_{2}),(k_{6},k_{3}),(k_{6},k_{4}),(k_{3},k_{4}),(k_{7},k_{5})\}.$$

We complete the proof noticing that  $H_{s}^{\rho_{1}} \oplus G^{\rho_{2}}=H_{t}^{\rho_{3}} \oplus G^{\rho_{2}}$. Thus every $\alpha$-galaxy is regular. 

\bbox

\section{Conclusions}

In that paper we showed new methods to prove the Erd\H{o}s-Hajnal Conjecture
for many classes of tournaments. In particular we showed that the Conjecture for the family of galaxies is s simple implication of the strong EH-property.
Strong EH-property is also useful to prove the Conjecture for a tournament $C_{5}$. Even though we did not show it, it can be proven using our techniques that every constellation satisfies the Conjecture. In other words, our techniques enable us to prove the Conjecture for all prime tournaments for which it was known so far. Moreover, we presented new infinite families of tournaments, containing infinitely many prime tournaments, with the EH-property.
It is interesting to note that the Conjecture for all tournaments on at most six vertices other than $H_{D}$ also follows from the strong EH-property.
The examples we gave in the previous section is just a small subset of all possible applications. By combining all new tournaments fow which we proved the Conjecture in this paper, using $\oplus$ operation and Theorem~\ref{extension_theorem}, we can prove the Conjecture for many other tournaments. 
We did not include other examples in the paper since to prove the Conjecture for them we use exactly the same techniques that turned out to work for families of tournaments we focused on here. The Reader probably realizes now how those techniques can be applied in other settings.
There are so many of them that a natural question arises. How can we characterize all of them ? Is there a compact description as it was in the case of galaxies or constellations ?
It is also worth to mention a striking difference between the way we understand the Conjecture in the directed and undirected setting. In the undirected one, despite many efforts, the Conjecture is not known for prime tournaments with more than $5$ vertices. This paper shows, that in the directed setting there are infinitely many prime tournaments that are not galaxies and even not constellations, but they satisfy the Conjecture.

%%%%%%%%%%%%%%%%%%
%Let $(v_{1},\eta_{1}) \in (v_{2},\eta_{2})$. Denote $v_{1} = %(v_2(i_{1}),...,v_{2}(i_{t}))$ for some $i_{1} < ... < i_{t}$.
%Let $f_{1}:V(H) \rightarrow \{1,2,...,k_{1}\}$, where: $k_{1}=\sum_{j} %\eta_{1}(j) + \sum_{j} [v_{1}(j)=0]$.
%We say that function $f_{2}:V(H) \rightarrow \{1,...,k_{2}\}$, where
%$k_{2} = \sum_{j} \eta_{2}(j) + \sum_{j}[v_{2}(j)=0]$ 
%%%%%%%%%%%%%%%%%%%%%%%%

\end{document}